\documentclass[final,10pt,3p]{elsarticle}

\usepackage[T1]{fontenc}
\usepackage[french,english]{babel}

\journal{Journal de Math\'ematiques Pures et Appliqu\'ees}

\usepackage{amsmath,amsthm}
\usepackage{amssymb}
\usepackage{amsfonts}
\usepackage{mathrsfs}

\newtheorem{theorem}{Theorem}[section]

\newtheorem{remark}{Remark}[section]

\numberwithin{equation}{section}

\newcommand{\erf}{{\rm erf}}
\newcommand{\ds}{\displaystyle}

\usepackage{color}
\usepackage{lineno,hyperref}
\modulolinenumbers[5]
\hypersetup{
	colorlinks = true,%
	citecolor = [rgb]{0.0,0.0,0.9},
	filecolor=black,%
	linkcolor = [rgb]{0.65,0.0,0.0},%
	anchorcolor = red,
	urlcolor= [rgb]{0.65,0.0,0.0}
}

\newenvironment{abstracts}
 {\global\setbox\absbox=\vbox\bgroup
    \hsize=\textwidth
    \linespread{1}\selectfont}
 {\vspace{-\bigskipamount}\egroup}
\renewenvironment{abstract}[1][]
 {\if\relax\detokenize{#1}\relax\else\selectlanguage{#1}\fi
  \noindent\textbf{\abstractname}\par\medskip\noindent\ignorespaces}
 {\par\bigskip}

 \newsavebox\extrainfobox

\begin{document}

\begin{frontmatter}

\title{Sharp uncertainty principles on Riemannian manifolds: the influence of curvature}

\author{Alexandru Krist\'aly}
 \address{Institute of Applied Mathematics,  \'Obuda
	University, 1034
	Budapest, Hungary}
\address{Department of Economics, Babe\c s-Bolyai University, 400591 Cluj-Napoca,
	Romania}
\ead{alex.kristaly@econ.ubbcluj.ro; alexandrukristaly@yahoo.com}

\begin{abstracts}
\begin{abstract}
  We present a  rigidity scenario for complete Riemannian manifolds
  supporting the Heisenberg-Pauli-Weyl uncertainty principle with the
  sharp constant in $\mathbb R^n$ (shortly, {\it sharp HPW
  	principle}). Our results deeply depend on the curvature of the Riemannian manifold which can be roughly formulated as follows: 
  \begin{itemize}
 \item[(a)] When $(M,g)$ has {\it non-positive sectional
  	curvature}, the sharp HPW principle holds on $(M,g)$. However,  {\it positive extremals exist} in the sharp HPW
  principle if and only if $(M,g)$ is
  isometric to $\mathbb R^n$,  $n={\rm dim}(M)$.

  \item[(b)] When $(M,g)$ has {\it non-negative Ricci curvature},  the sharp HPW principle holds on $(M,g)$ if and only if $(M,g)$ is
  isometric to $\mathbb R^n$.
\end{itemize}
  Since the sharp HPW principle and the Hardy-Poincar\'e inequality
  are endpoints of the Caffarelli-Kohn-Nirenberg interpolation
  inequality, we establish further quantitative results for the latter
  inequalities in terms of the curvature on Cartan-Hadamard manifolds. 
\end{abstract}

\begin{abstract}[french]
	Nous pr\'esentons un sc\'enario de rigidit\'e pour les vari\'et\'es riemanniennes compl\`etes soutenant le principe d'incertitude d'Heisenberg-Pauli-Weyl avec la constante optimale en $\mathbb R^n$ (bri\`evement, {\it le principle d'HPW}). Nos r\'esultats d\'ependent profond\'ement de la courbure de la vari\'et\'e riemannienne et ils peuvent \^etre formul\'es comme suit:
	
	\begin{itemize}
		\item[(a)] Lorsque $(M, g)$ a {\it courbure sectionnelle non positive}, le principe d'HPW a lieu sur $(M,g)$. N\'eanmoins, {\it des fonctions extr\'emales positives existent} dans le principe d'HPW si et seulement si $(M, g)$ est isom\'etrique \`a $\mathbb R^n$, $n = {\rm dim} (M) $.
		\item[(b)] Lorsque $(M,g)$ a {\it courbure de Ricci non n\'egative}, le principe d'HPW a lieu sur $(M,g)$ si et seulement si $(M,g)$ est isom\'etrique \`a $\mathbb R^n$.
	\end{itemize}
	Comme le principe d'HPW et l'in\'egalit\'e Hardy-Poincar\'e sont des cas extr\^emes de l'in\'egalit\'e d'interpolation de Caffarelli-Kohn-Nirenberg, nous \'etablissons des r\'esultats quantitatifs pour les derni\`eres in\'egalit\'es en terme de la courbure sur les vari\'et\'es de Cartan-Hadamard.
\end{abstract}
\end{abstracts}

\begin{keyword}
Heisenberg-Pauli-Weyl uncertainty principle; Riemannian
	manifold; sharp constant; curvature
\MSC[2010]{Primary 53C21; Secondary 58J60}
\end{keyword}

\tnotetext[t1]{
Dedicated to professor Zolt\'an M. Balogh on the occasion of his 50th birthday
}
\tnotetext[t1]{
	Research supported by the Project CNFIS-FDI-2016-0056, STAR-UBB Fellowship.
}

\end{frontmatter}


 \section{Introduction and main results}

The Heisenberg uncertainty principle in quantum mechanics states
that the position and momentum of a given particle cannot be
accurately determined simultaneously, see \cite{Heisenberg}. The
rigorous mathematical formulation of this principle is attributed to
Pauli and Weyl \cite{PauliWeyl}, stating that the function itself
and its Fourier transform cannot be sharply localized at the same
time. In terms of PDEs, the Heisenberg-Pauli-Weyl uncertainty
principle in the Euclidean setting is described by the inequality
{\small \begin{equation}\label{euklideszi-HPW}
	\left(\int_{\mathbb R^n}|\nabla u(x)|^2{\text d}x\right) \left(\int_{\mathbb
		R^n}|x|^2u(x)^2{\text d}x\right)\geq \frac{n^2}{4}\left(\int_{\mathbb
		R^n}u(x)^2{\text d}x\right)^2, \ \forall  u\in C_0^\infty(\mathbb R^n).
	\end{equation}}
It is well known that the constant $\frac{n^2}{4}$ is sharp and the
extremals are given (up to a constant) by the family of  Gaussian
functions $u_\lambda(x)=e^{-\lambda |x|^2}$, $\lambda>0.$

Since its initial formulation, the Heisenberg-Pauli-Weyl principle
is deserving continuously a deep source of inspiration in different
areas of Physics and Mathematics. Without the sake of completeness, the
Heisenberg-Pauli-Weyl principle has been studied in various
contexts, see  Ciatti, Ricci and Sundari \cite{Ciatti} (for positive
self-adjoint operators on measure spaces), Fefferman \cite{Fef},
Folland and Sitaram \cite{Folland}, and Nahmod \cite{Nahmod}
(locating eigenvalues for selfadjoint differential operators via SAK
principle), Andersen \cite{Andersen, Andersen-2}, Erb \cite{Erd, Erb-Phd}  and Kombe and \"Ozaydin
\cite{KO-TAMS-2009, KO-TAMS-2013} (sharp uncertainty principle on
compact/noncompact Riemannian manifolds), Okoudjou, Saloff-Coste and
Teplyaev \cite{Ok} (for fractals, graphs and metric measure spaces),
and references therein.

The  purpose of our paper is to describe a complete scenario
concerning the {sharp} Heisen\-berg-Pauli-Weyl uncertainty
principle on complete Riemannian manifolds. 
Hereafter, in order to
avoid confusions,  the {\it sharpness} is understood in the sense
that the Heisenberg-Pauli-Weyl principle holds on a Riemannian
manifold $(M,g)$ with the same constant $\frac{n^2}{4}$ as in the
Euclidean space $\mathbb R^n$.

To be more precise,  let $(M,g)$ be an $n(\geq 2)-$dimensional
complete Riemannian manifold, ${\text
	d}V_g$ its canonical volume element, and $d_{x_0}(x)=d(x_0,x)$ be the distance function from a
point $x_0\in M$. For $x_0\in M$ fixed, we consider the {\it
	Hei\-senberg-Pauli-Weyl principle} on $(M,g)$ of the
form: for all $u\in C_0^\infty(M)$,
{$$
	\left(\int_{M}|\nabla_gu|^2{\text d}V_g\right) \left(\int_{M}d_{x_0}^2u^2{\text d}V_g\right) \geq \frac{n^2}{4}\left(\int_{M}u^2{\text d}V_g\right)^2. \eqno{{\bf ({HPW})}_{x_0}}
	$$}

Our first result can be stated as follows:

\begin{theorem}\label{theorem-uncertainty-intro-1}{\rm [Non-positively curved case]} Let $(M,g)$ be
	an
	$n-$di\-men\-sional Car\-tan-Hadamard manifold $($simply connected, complete Riemannian manifold with non-positive sectional curvature$)$.
	\begin{itemize}
		\item[{\rm (i)}] {\rm [Sharpness]} The Heisenberg-Pauli-Weyl principle  ${\bf
			({HPW})}_{x_0}$ holds for every $x_0\in M;$ moreover,
		$\frac{n^2}{4}$ is sharp, i.e.,{
			$$\frac{n^2}{4}=\inf_{u\in C_0^\infty(M)\setminus \{0\}}{ \frac{\left(\displaystyle\int_{M}|\nabla_gu|^2{\rm d}V_g\right)
					\left(\displaystyle\int_{M}d_{x_0}^2u^2{\rm
						d}V_g\right)}{\left(\displaystyle\int_{M}u^2{\rm
						d}V_g\right)^2}}.
			$$
		}
		\item[{\rm (ii)}] {\rm [Extremals]}  The following
		statements are equivalent:
		\begin{itemize}
			\item[\rm (a)]  $\frac{n^2}{4}$ is achieved by a positive extremal  in ${\bf ({HPW})}_{x_0}$ for {\rm some} $x_0\in M;$
			\item[\rm (b)]   $\frac{n^2}{4}$ is achieved by a positive extremal  in ${\bf ({HPW})}_{x_0}$ for {\rm every} $x_0\in M;$
			\item[\rm (c)] $(M,g)$ is isometric to $\mathbb R^n$.
		\end{itemize}
	\end{itemize}
\end{theorem}

Some remarks are in order concerning Theorem \ref{theorem-uncertainty-intro-1}. 
\begin{remark}\rm
	(a) ${\bf ({HPW})}_{x_0}$ is a consequence of a quantitative/weighted  Heisenberg-Pauli-Weyl  principle stated below in Theorem  \ref{quantitative-HPW}. Note that similar weighted Heisenberg-Pauli-Weyl type principles have been investigated on Riemannian manifolds diffeomorphic to $\mathbb R^n$ (thus, in particular, on Cartan-Hadamard manifolds). Indeed, by using an operator theoretic approach, Erb \cite[Theorem 2.54]{Erb-Phd} stated weighted Heisenberg-Pauli-Weyl principles where the weights are in terms of volume distortion coefficients involving information on the curvature of the manifold. By using Bishop-Gromov comparison arguments, Corollary 2.68 of Erb \cite{Erb-Phd} can be seen as ${\bf ({HPW})}_{x_0}$.   
	Note however that the sharpness and the characterization of extremals in ${\bf ({HPW})}_{x_0}$ are not explicitly investigated in \cite{Erb-Phd}.

	(b) One could expect finer results for  ${\bf
		({HPW})}_{x_0}$ whenever the Riemannian manifold is the model hyperbolic space. Andersen \cite{Andersen, Andersen-2} proved that hyperbolic Gaussians are candidates for extremal functions in  Heisenberg-Pauli-Weyl principles within the hyperbolic setting. Recently,  Kombe and \"Ozaydin \cite[Theorem 4.2]{KO-TAMS-2013} claimed that the hyperbolic Gaussian
	function $u(x)=e^{-\alpha d(x)^2}$ (where $\alpha>0$ is a root of a
	highly nonlinear equation) is an extremal function in the
	Heisenberg-Pauli-Weyl principle  ${\bf ({HPW})}_0$ on the hyperbolic
	space $\mathbb H^n$; hereafter,   $d(x)=d_{\mathbb H^n}(0,x)$ denotes the hyperbolic distance between $0$ and $x$ in the Poincar\'e ball model. 
	According to Theorem \ref{theorem-uncertainty-intro-1},  the scenario
	described in \cite{KO-TAMS-2013} cannot occur; moreover, two further
	independent arguments are presented in \S \ref{subsec-4-3} which
	confirm the fact that in the hyperbolic setting the expected
	Gaussian function $u(x)=e^{-\alpha d(x)^2}$ {\it cannot} be extremal
	in ${\bf ({HPW})}_0$ for {\it any} $\alpha>0.$ More precisely, the hyperbolic Gaussians are extremals for a quantitative Heisenberg-Pauli-Weyl principle rather than for ${\bf ({HPW})}_0$, as we shall explain in the sequel, see (\ref{uncer-qUANT1}). 
	
	(c) Being within the context of Cartan-Hadamard manifolds, the sharpness of  Sobolev-type inequalities  usually requires the validity of the longstanding Cartan-Ha\-da\-mard conjecture, i.e., the sharp isoperimetric inequality (which is valid in $2,$ $3$ and $4-$dimensional Cartan-Hadamard manifolds), see e.g. Hebey \cite[Section 8.2]{Hebey}. We notice that such a hypothesis is not needed in Theorem \ref{theorem-uncertainty-intro-1}. 
\end{remark}

In the  non-negatively curved case the situation is even more rigid than in Theorem \ref{theorem-uncertainty-intro-1}:

\begin{theorem}\label{theorem-uncertainty-intro-2}{\rm [Non-negatively curved case]} Let $(M,g)$ be a  complete,
	$n-$di\-men\-sional Riemannian manifold with non-negative Ricci curvature. The following
	statements are equivalent:
	\begin{itemize}
		\item[\rm (a)]   ${\bf
			({HPW})}_{x_0}$ holds for {\rm some} $x_0\in M;$
		\item[\rm (b)]   ${\bf
			({HPW})}_{x_0}$ holds for {\rm every} $x_0\in M;$
		\item[\rm (c)] $(M,g)$ is isometric to $\mathbb R^n$.
	\end{itemize}
\end{theorem}

\begin{remark}\rm
	Theorem \ref{theorem-uncertainty-intro-2} can be included into the
	{\it best constant program} initiated by Aubin \cite{Aubin},  and
	studied by Ledoux \cite{Ledoux}, Cheeger and Colding \cite{CC},
	Druet, Hebey and Vaugon \cite{Druetetal},  do Carmo and Xia
	\cite{doCarmo-Xia}, Minerbe \cite{Minerbe}, Li and Wang
	\cite{Li-Wang}, etc. Indeed, in the aforementioned papers, the
	authors established that complete Riemannian manifolds with
	non-negative Ricci curvature supporting some Sobolev-type
	inequalities should be close to Euclidean spaces whenever the constant is sufficiently close to the sharp Euclidean Sobolev constant. The reader may
	consult Hebey \cite{Hebey} for a thoroughgoing presentation of this
	subject.
\end{remark}

In the sequel, we shall present some closely related results to the
sharp Heisenberg-Pauli-Weyl principle on Riemannian manifolds which
are of independent interests.

Let $p,q\in \mathbb R$ and $n\in \mathbb N$
be such that
\begin{equation}\label{p_q_feltetelek}
0<q<2<p\ {\rm and}\ 2<n<\frac{2(p-q)}{p-2}.
\end{equation}
For a fixed $x_0\in M$, we consider the  {\it Caffarelli-Kohn-Nirenberg interpolation
	inequality} on $(M,g)$:  for all $u\in C_0^\infty(M)$, {
	$$
	\left(\int_{M}|\nabla_gu|^2{\text d}V_g\right) \left(\int_{M}\frac{|u|^{2p-2}}{d_{x_0}^{2q-2}}{\text d}V_g\right) \geq \frac{(n-q)^2}{p^2}\left(\int_{M}
	\frac{|u|^p}{d_{x_0}^q} {\text d}V_g\right)^2.\eqno{{\bf (CKN)}_{x_0}}
	$$}

An {endpoint} of ${\bf (CKN)}_{x_0}$ is precisely the 
Heisenberg-Pauli-Weyl principle ${\bf ({HPW})}_{x_0}$ whenever $p\to
2$ and $q\to 0$. As a part of the {\it best constant program,} Xia
\cite{Xia} proved that if $(M,g)$ is a complete, $n-$di\-men\-sional
Riemannian manifold with non-negative Ricci curvature,  then $(M,g)$
supports ${\bf (CKN)}_{x_0}$ for some $x_0\in M$ if and only if
$(M,g)$ is isometric to $\mathbb R^n$. In the Euclidean setting, Xia
\cite{Xia} also proved the sharpness of $\frac{(n-q)^2}{p^2} $ in
${\bf (CKN)}_{x_0}$ and the existence of a class of extremals
\begin{equation}\label{caff-extremal}
u_\lambda(x)={\left(\lambda+|x-x_0|^{2-q}\right)^{\frac{1}{2-p}}},\
\lambda>0.
\end{equation}
The reader may also consult  Krist\'aly and Ohta
\cite{Kristaly-Ohta} for a study of Caffarelli-Kohn-Nirenberg
inequalities on 'positively curved' metric measure spaces.

The non-positively curved counterpart of Xia's result, similar to
Theorem \ref{theorem-uncertainty-intro-1}, can be stated as follows:

\begin{theorem}\label{theorem-interpolation-intro}
	Let  $p,q\in \mathbb R$ and $n\in \mathbb N$
	be such that {\rm (\ref{p_q_feltetelek})} holds and let $(M,g)$ be an
	$n-$di\-men\-sional Car\-tan-Hadamard manifold.
	\begin{itemize}
		\item[{\rm (i)}] {\rm [Sharpness]} The Caffarelli-Kohn-Nirenberg interpolation
		inequality  ${\bf ({CKN})}_{x_0}$ holds for every $x_0\in M$ and the
		constant  $\frac{(n-q)^2}{p^2}$ is sharp, i.e., {\small
			$$\frac{(n-q)^2}{p^2}=\inf_{u\in C_0^\infty(M)\setminus \{0\}}{ \frac{\left(\displaystyle\int_{M}|\nabla_gu|^2{\rm d}V_g\right)
					\left(\displaystyle\int_{M}\frac{|u|^{2p-2}}{d_{x_0}^{2q-2}}{\rm
						d}V_g\right)}{\left(\displaystyle\int_{M}\frac{|u|^p}{d_{x_0}^q}{\rm
						d}V_g\right)^2}}.$$}
		\item[{\rm (ii)}] {\rm [Extremals]}  The following
		statements are equivalent:
		\begin{itemize}
			\item[\rm (a)]  {\small $\frac{(n-q)^2}{p^2}$} is achieved by a positive extremal  in ${\bf ({CKN})}_{x_0}$ for {\rm some} {\small $x_0\in
				M;$}
			\item[\rm (b)]   {\small $\frac{(n-q)^2}{p^2}$} is achieved by a positive extremal  in ${\bf ({CKN})}_{x_0}$ for {\rm every} {\small $x_0\in
				M;$}
			\item[\rm (c)] $(M,g)$ is isometric to $\mathbb R^n$.
		\end{itemize}
	\end{itemize}
\end{theorem}

The other endpoint of ${\bf (CKN)}_{x_0}$, whenever $p\to 2$ and
$q\to 2$, is the famous {\it Hardy-Poincar\'e inequality} on
$(M,g)$: for all $u\in C_0^\infty(M)$, {
	$$
	\int_{M}|\nabla_gu|^2{\text d}V_g  \geq \frac{(n-2)^2}{4}\int_{M}
	\frac{u^2}{d_{x_0}^2} {\text d}V_g.\eqno{{\bf (HP)}_{x_0}}
	$$}
In the Euclidean setting it is well known that $\frac{(n-2)^2}{4}$
is sharp, but there are no extremal functions. The lack of extremals
motivated various improvements of the Hardy-Poincar\'e inequality;
see e.g. Adimurthi, Chaudhuri and Ramaswamy \cite{AR}, Barbatis,
Filippas and Tertikas \cite{Barbatis},   Brezis and V\'azquez
\cite{BV}, Filippas and Tertikas \cite{Fil-Ter}, Ghoussoub and
Moradifam \cite{Gho-Morad, Gho-Morad-2}, Wang and Willem \cite{WW},
etc.

In the last few years, the Hardy-Poincar\'e inequality has been also
studied on complete, non-compact Riemannian manifolds, where the
influence of geometry played a key role; see e.g. Berchio,  D'Ambrosio,  Ganguly and Grillo \cite{Berchio}, Berchio,   Ganguly and Grillo \cite{Berchio-2}, Carron 
\cite{Carron-JMPA}, D'Ambrosio and Dipierro \cite{D-D-olaszok},
Kombe and \"Ozaydin \cite{KO-TAMS-2009,KO-TAMS-2013}, Yang, Su and
Kong \cite{YSK-Kinaiak}, and references therein. 

Our aim is to provide a new type of improved Hardy-Poincar\'e inequality
which shows that {\it more curvature implies more powerful
	improvements}:

\begin{theorem}\label{theorem-Hardy-intro} {\rm [Improved Hardy-Poincar\'e inequality via curvature]}
	Let $(M,g)$ be an  $n-$di\-men\-sional Cartan-Hadamard manifold such
	that the sectional curvature is bounded from above by $c\leq 0.$
	Then for every $x_0\in M$ and $u\in C_0^\infty(M),$ we have 
	$$
	\int_{M}|\nabla_gu|^2{\rm d}V_g  \geq  \frac{(n-2)^2}{4}\int_{M}
	\frac{u^2}{d_{x_0}^2} {\rm d}V_g  +\frac{3{|c|}(n-1)(n-2)}{2}\int_M \frac{u^2}{\pi^2+|c|d_{x_0}^2}{\rm{d}}V_g.
	$$ In addition,  the constant $\frac{(n-2)^2}{4}$ is sharp $($independently by the second term on the RHS$).$
\end{theorem}

\begin{remark}\rm  It seems similar rigidity results for the Hardy-Poincar\'e inequalities as in the Theorem \ref{theorem-uncertainty-intro-2} cannot be established on non-negatively curved spaces.
	In the proof of Theorem \ref{theorem-uncertainty-intro-2} the
	existence of extremals in the Euclidean case is crucial which fails
	in the case of Hardy-Poincar\'e inequalities.

\end{remark}

{\it Plan of the paper.} In Section \ref{sect-2} we first recall the
notions and results from Riemannian geometry which are used
throughout the proofs. In Section \ref{sect-3} we first deal with
the generic Heisenberg-Pauli-Weyl principle by proving Theorems
\ref{theorem-uncertainty-intro-1}\&\ref{theorem-uncertainty-intro-2}, and then we consider this
principle on hyperbolic spaces (w.r.t. the paper
\cite{KO-TAMS-2013}). In Section \ref{sect-4} we study related
inequalities to the Heisenberg-Pauli-Weyl principle on
Cartan-Hadamard manifolds (i.e., Caffarelli-Kohn-Nirenberg
interpolation inequality and Hardy-Poincar\'e inequality).

\section{Preliminaries}\label{sect-2}

Let $(M,g)$ be an $n-$dimensional complete Riemannian manifold, and
$d:M\times M\to [0,\infty)$ be the metric function associated to the
Riemannian metric $g$. Let $B(x,\rho)=\{y\in M:d(x,y)<\rho\}$ be the
open metric ball with center $x\in M$ and radius $\rho>0.$ If
${\text d}V_g$ is the canonical volume
element on $(M,g)$, the volume of a bounded open set $S\subset M$ is
$${\rm Vol}_g(S)=\displaystyle\int_S {\text d}V_g={\rm Haus}_{d}(S),$$ where ${\rm
	Haus}_{d}(S)$ is the Hausdorff measure of $S$ with respect to the
metric function $d$.  In general, one has for every $x\in M$ that
\begin{equation}\label{volume-comp-nullaban}
\lim_{\rho\to 0^+}\frac{{\rm Vol}_g(B(x,\rho))}{\omega_n
	\rho^n}=1,
\end{equation}
where $\omega_n$ is the volume of the standard $n-$dimensional
Euclidean unit ball.

Let $u:M\to \mathbb R$ be of class $C^1.$ If $(x^i)$ is the local
coordinate system on a coordinate neighborhood of $x\in M$, and the
local components of the differential of $u$ are denoted
$u_i=\frac{\partial u}{\partial x_i}$, then the local components of
the gradient  $\nabla_g u$ are $u^i=g^{ij}u_j$. Here, $g^{ij}$ are
the local components of $g^{-1}=(g_{ij})^{-1}$.

The Laplace-Beltrami operator is given by $\Delta_g u={\rm
	div}(\nabla_g u)$ whose expression in a local chart of associated
coordinates $(x^i)$ is $$\Delta_g u=g^{ij}\left(\frac{\partial^2
	u}{\partial x_i\partial x_j}-\Gamma_{ij}^k\frac{\partial u}{\partial
	x_k}\right),$$ where $\Gamma_{ij}^k$ are the coefficients of the
Levi-Civita connection. 

If $u,v:M\to \mathbb R$ are of class $C^2$, one has the following
integration by parts formula
$$\int_M v \Delta_g u {\text d}V_g=-\int_M \langle \nabla_g v,\nabla_g u\rangle {\text d}V_g,$$
where $\langle\cdot,\cdot \rangle$ denotes the scalar product
associated with the Riemannian metric $g$ for $1-$forms. For
simplicity, we shall use the notation
$|\alpha|=\sqrt{\langle\alpha,\alpha \rangle}$ for any $1-$form.

A Riemannian manifold $(M,g)$ is called Cartan-Hadamard if it is
complete, simply connected and with non-positive sectional
curvature.

For every $c\leq 0$ we consider the function ${\bf
	ct}_c:(0,\infty)\to \mathbb R$ defined by
$${\bf ct}_c(\rho)=\left\{
\begin{array}{lll}
\frac{1}{\rho}
& \hbox{if} &  c=0, \\
\sqrt{|c|}\coth(\sqrt{|c|}\rho) & \hbox{if} & c<0.
\end{array}\right.$$
For further use, let ${\bf D}_c:[0,\infty)\to \mathbb R$ defined by
$${\bf D}_c(\rho)=\left\{
\begin{array}{lll}
0
& \hbox{if} &  \rho=0, \\
\rho {\bf ct}_c(\rho)-1 & \hbox{if} & \rho>0.
\end{array}\right.$$
It is clear that ${\bf D}_c\geq 0.$

Hereafter, $d_{x_0}(x)=d(x_0,x)$ denotes the distance function from
a given point $x_0\in M$.

\begin{theorem}\label{comparison-laplace}{\rm [Laplacian comparison; see \cite[Theorem 5.1]{Wu-Xin}]} Let $(M,g)$ be an $n-$di\-men\-sional Cartan-Hadamard manifold
	such that the sectional curvature is bounded from above by $c\leq
	0$, and let $x_0\in M$ be fixed. Then we have $($in distributional
	sense$)$ that
	$${\Delta}_gd_{x_0}\geq (n-1){\bf ct}_c(d_{x_0}).$$
\end{theorem}

In the proof of our results Bishop-Gromov-type volume comparison
principles  play a crucial role. Here we adapt from the Finsler
version the following form (see Shen \cite{Shen-volume}, Wu and Xin
\cite[Theorems 6.1 \& 6.3]{Wu-Xin} and Zhao and Shen
\cite{Zhao-Shen}):

\begin{theorem}\label{comparison-volume}{\rm [Volume comparison]}
	Let $(M,g)$ be a complete, $n-$di\-men\-sional Riemannian manifold.
	Then the following statements hold.
	\begin{itemize}
		\item[{\rm (a)}] If $(M,g)$ is a Cartan-Hadamard manifolds,
		the function
		$\rho\mapsto \frac{{\rm Vol}_g(B(x,\rho))}{\rho^n}$ is
		non-decreasing, $\rho>0$. 
		In particular, from {\rm (\ref{volume-comp-nullaban})} we have
		\begin{equation}\label{volume-comp-altalanos-0}
		{{\rm Vol}_g(B(x,\rho))}\geq \omega_n \rho^n\ {for\ all}\ x\in M\
		{and}\ \rho>0.
		\end{equation}
		If equality holds in {\rm (\ref{volume-comp-altalanos-0})}, then the
		sectional curvature is identically zero.
		\item[{\rm (b)}] If $(M,g)$ has non-negative Ricci curvature,
		the function
		$\rho\mapsto \frac{{\rm Vol}_g(B(x,\rho))}{\rho^n}$ is
		non-increasing, $\rho>0$. In particular, from {\rm
			(\ref{volume-comp-nullaban})} we have
		\begin{equation}\label{volume-comp-altalanos-2}
		{{\rm Vol}_g(B(x,\rho))}\leq \omega_n \rho^n\ {for\ all}\ x\in M\
		{and}\ \rho>0.
		\end{equation}
		If equality holds in {\rm (\ref{volume-comp-altalanos-2})}, then the
		sectional curvature is identically zero.
	\end{itemize}
\end{theorem}

\vspace{0.5cm}

\section{Heisenberg-Pauli-Weyl principle on Riemannian
	manifolds}\label{sect-3}

\subsection{Non-positively curved case:  proof of Theorem \ref{theorem-uncertainty-intro-1}}

First, we present a quantitative version of the
Heisenberg-Pauli-Weyl principle.

\begin{theorem}{\rm [Quantitative
		Heisenberg-Pauli-Weyl principle]}\label{quantitative-HPW} Let
	$(M,g)$ be an $n-$di\-men\-sional Cartan-Hadamard manifold such that
	the sectional curvature is bounded from above by $c\leq 0.$ Then for
	all $x_0\in M$ and $u\in C_0^\infty(M)$, we have
	$$
	\left(\int_{M}|\nabla_g u|^2{\rm d}V_g\right) \left(\int_Md_{x_0}^2u^2{\rm d}V_g\right)\geq
	\frac{n^2}{4}\left(\int_M\left(1+\frac{n-1}{n}{\bf
		D}_c(d_{x_0})\right)u^2{\rm d}V_g\right)^2.
	$$
\end{theorem}

{\it Proof.} Let  $x_0\in M$ and  $u\in C_0^\infty(M)$ be fixed
arbitrarily. According to  Theorem \ref{comparison-laplace}, one has
\begin{eqnarray}
\label{laplace-egyenlotlenseg}
\int_M  {\Delta}_g(d_{x_0}^2)u^2{\text d}V_g &=& 2\int_M (1+d_{x_0}{\Delta}_gd_{x_0})u^2
{\text d}V_g \nonumber \\
&\geq& 2 \int_M  (1+(n-1)d_{x_0}{\bf ct}_c(d_{x_0}))u^2{\text
	d}V_g \nonumber\\
&=&2n \int_M  \left(1+\frac{n-1}{n}{\bf D}_{c}(d_{x_0})\right)u^2{\text
	d}V_g.
\end{eqnarray}
An integration by parts yields
\begin{eqnarray*}
	\int_M {\Delta}_g(d_{x_0}^2)u^2{\text d}V_g &=& -\int_M \langle \nabla_g (u^2),{\nabla_g}(d_{x_0}^2)\rangle{\text d}V_g \\
	&=& -4\int_M
	ud_{x_0}\langle \nabla_g u,{\nabla_g}d_{x_0}\rangle{\text d}V_g.
\end{eqnarray*}
By the eikonal equation $|{\nabla}_gd_{x_0}|=1$ a.e. on $M$, one has that $|\langle \nabla_g
u,{\nabla_g}d_{x_0}\rangle|\leq |\nabla_g u|.$ Thus, by Schwartz
inequality one gets { $$\left(\int_M
	ud_{x_0}\langle \nabla_g u,{\nabla_g}d_{x_0}\rangle{\text d}V_g\right)^2\leq
	\left(\int_{M}d_{x_0}^2u^2{\text d}V_g\right)\left(\int_{M}|\nabla_g u|^2{\text d}V_g\right).
	$$}
The latter relation coupled with (\ref{laplace-egyenlotlenseg})
yields the quantitative Heisenberg-Pauli-Weyl principle, which
concludes the proof.
\hfill $\square$\\

{\it Proof of Theorem \ref{theorem-uncertainty-intro-1}.} (i) Let
$x_0\in M$ be fixed. Since ${\bf D}_c\geq 0,$ due to Theorem
\ref{quantitative-HPW}, the Heisenberg-Pauli-Weyl principle ${\bf
	({HPW})}_{x_0}$ holds.

We shall prove that the constant $\frac{n^2}{4}$ is optimal in ${\bf
	({HPW})}_{x_0}$, following Aubin's argument \cite{Aubin}; see also
Hebey \cite{Hebey}. Let
\begin{equation}\label{c-hpw}
\textsf{C}_\textsf{HPW}=\inf_{u\in C_0^\infty(M)\setminus \{0\}}{
	\frac{\left(\displaystyle\int_{M}|\nabla_gu|^2{\rm d}V_g\right)
		\left(\displaystyle\int_{M}d_{x_0}^2u^2{\rm
			d}V_g\right)}{\left(\displaystyle\int_{M}u^2{\rm d}V_g\right)^2}}.
\end{equation}
Since ${\bf ({HPW})}_{x_0}$ holds, then $\textsf{C}_\textsf{HPW}\geq
\frac{n^2}{4}$. Assume  that
$\textsf{C}_\textsf{HPW}> \frac{n^2}{4}$. By (\ref{c-hpw}), one has
\begin{equation}\label{c-hpw-uj}
\left(\int_{M}|\nabla_gu|^2{\text d}V_g\right) \left(\int_{M}d_{x_0}^2u^2{\text d}V_g\right) \geq \textsf{C}_\textsf{HPW}\left(\int_{M}u^2{\text d}V_g\right)^2,\ \forall u\in C_0^\infty(M).
\end{equation}
For every $\varepsilon>0$, there exists a local chart
$(\Omega,\phi)$ of $M$ at the point $x_0$ and a number $\delta>0$
such that $\phi(\Omega)=B_e(0,\delta)$ and the components $g_{ij}$
of the metric $g$ satisfy
\begin{equation}\label{two-sided}
(1-\varepsilon)\delta_{ij}\leq g_{ij} \leq (1+\varepsilon)\delta_{ij}
\end{equation}
in the sense of bilinear forms. Here, $B_e(0,\delta)$ is the
$n-$dimensional Euclidean ball of center $0$ and radius $\delta>0$.

According to (\ref{c-hpw-uj}) and to the two-sided metric estimate
(\ref{two-sided}), for $\varepsilon>0$ small enough, there exists
$\tilde \delta>0$ and $\textsf{C}_\textsf{HPW}'>\frac{n^2}{4}$ such
that for every $\delta\in (0,\tilde \delta)$ and $w\in
C_0^\infty(B_e(0,\delta))$,
\begin{equation}\label{c-hpw-uj-meg}
\left(\int_{B_e(0,\delta)}|\nabla w|^2{\text d}x\right) \left(\int_{B_e(0,\delta)}|x|^2w^2{\text d}x\right) \geq \textsf{C}_\textsf{HPW}'\left(\int_{B_e(0,\delta)}w^2{\text d}x\right)^2.
\end{equation}
Let $u\in C_0^\infty(\mathbb R^n)$ be fixed arbitrarily and set
$w_\lambda(x)=u(\lambda x)$, $\lambda>0.$ It is clear that
$w_\lambda\in C_0^\infty(B_e(0,\delta))$ for enough large
$\lambda>0.$ Inserting $w_\lambda$ into (\ref{c-hpw-uj-meg}), and
having the scaling properties
$$\int_{B_e(0,\delta)}|\nabla w_\lambda|^2{\text d}x=\lambda^{2-n}\int_{\mathbb R^n}|\nabla u|^2{\text d}x,\ \int_{B_e(0,\delta)}|x|^2w_\lambda^2{\text d}x=\lambda^{-2-n}\int_{\mathbb R^n}|x|^2u^2{\text d}x,$$
and
$$\int_{B_e(0,\delta)}w_\lambda^2{\text d}x=\lambda^{-n}\int_{\mathbb R^n}u^2{\text d}x,$$
it follows that
$$ \left(\int_{\mathbb R^n}|\nabla u|^2{\text d}x\right) \left(\int_{\mathbb R^n}|x|^2u^2{\text d}x\right) \geq \textsf{C}_\textsf{HPW}'\left(\int_{\mathbb R^n}u^2{\text d}x\right)^2.$$
In particular, in the latter relation we may substitute the Gaussian function
$u(x)=e^{-|x|^2}$, obtaining that $\frac{n^2}{4}\geq
\textsf{C}_\textsf{HPW}',$ a contradiction. Consequently, $\textsf{C}_\textsf{HPW}=\frac{n^2}{4}$.

(ii) First, if $(M,g)$ is isometric to $\mathbb R^n,$ the sharp Heisenberg-Pauli-Weyl
principle ${\bf ({HPW})}_{x_0}$ can be equivalently transformed into
(\ref{euklideszi-HPW}) for which the Gaussian functions
$u_\lambda(x)=e^{-\lambda|x|^2},$ $\lambda>0$, are extremal
functions. Thus, the implications
(c)$\Rightarrow$(b)$\Rightarrow$(a) hold true.

We now prove (a)$\Rightarrow$(c). Let $u_0>0$ be an extremal
function in the sharp Heisenberg-Pauli-Weyl principle ${\bf
	({HPW})}_{x_0}$ for some $x_0\in M.$ In particular, in the estimates
in Theorem \ref{quantitative-HPW} we should have equalities; thus,
by (\ref{laplace-egyenlotlenseg}) one has ${\bf D}_c\equiv 0$ (i.e.,
we necessarily have $c=0$, so the sectional curvature of $(M,g)$
cannot be bounded above by a fixed negative number), and
\begin{equation}\label{delta-egyenloseg}
{\Delta}_g(d_{x_0}^2) = 2n.
\end{equation}

Let us fix $\rho>0$ arbitrarily. Note that the unit outward pointing
normal vector to the sphere $S(x_0,\rho)=\partial B(x_0,\rho)=\{x\in
M:d({x_0},x)=\rho\}$ is given by ${\bf n}={\nabla}_g d_{x_0}$. Let
us denote by  ${\text d}\varsigma_g$ the volume form on
$S(x_0,\rho)$ induced from ${\text d}V_g$. By applying Stokes'
formula 
and the fact that $\langle{\bf n},{\bf n}\rangle=1$ we have
\begin{eqnarray*}
	2n {\rm Vol}_g(B(x_0,\rho)) &=& \int_{B(x_0,\rho)}{\Delta}_g(d_{x_0}^2){\text d}V_g=\int_{B(x_0,\rho)}{\rm div}({\nabla}_g(d_{x_0}^2)){\text d}V_g \\
	&{=}& \int_{S(x_0,\rho)}  \langle{\bf n},{\nabla}_g(d_{x_0}^2)\rangle{\text d}\varsigma_g=2\int_{S(x_0,\rho)} d_{x_0} \langle{\bf n},{\nabla}_g d_{x_0}\rangle{\text d}\varsigma_g \\
	&=&2\rho\int_{S(x_0,\rho)} \langle{\bf n},{\bf n}\rangle{\text d}\varsigma_g=2\rho\int_{S(x_0,\rho)} {\text
		d}\varsigma_g\\&=&2\rho{\rm A}_g({S(x_0,\rho)}),
\end{eqnarray*}
where
$${\rm
	A}_g({S(x_0,\rho)})=\lim_{\varepsilon\to 0^+}\frac{{\rm
		Vol}_g(B(x_0,\rho+\varepsilon))-{\rm
		Vol}_g(B(x_0,\rho))}{\varepsilon}:=\frac{{\text d}}{{\text
		d}\rho}{\rm Vol}_g(B(x_0,\rho))$$ is the surface area of
$S(x_0,\rho)$. Thus, the above relations imply that
$$\frac{\frac{{\text d}}{{\text d}\rho}{\rm Vol}_g(B(x_0,\rho))}{{\rm
		Vol}_g(B(x_0,\rho))}=\frac{n}{\rho}.$$ By integrating this
expression and due to relation (\ref{volume-comp-nullaban}), we
conclude that
\begin{equation}\label{ideiskell}
{\rm Vol}_g(B(x_0,\rho))=\omega_n \rho^n\ {\rm for\ all}\ \rho>0.
\end{equation}

Let $x\in M$ and $\rho>0$ be arbitrarily fixed. Since $(M,g)$ is of
Cartan-Hadamard type, by the volume comparison (see Theorem
\ref{comparison-volume}(a)),  the function $r\mapsto \frac{{\rm
		Vol}_g(B(x,r))}{r^n}$ is non-decreasing on $(0,\infty)$. Therefore,
one has
\begin{eqnarray*}
	\omega_n &\leq&\frac{{\rm Vol}_g(B(x,\rho))}{\rho^n}\ \ \ \ \ \ \ \
	\ \ \
	\ \ \ \ \ \ \ \ \ \ \ \ \ \ \ \ \ \
	\ \ \ \ \ \ \ \ \ \ \ \ \ \ \ \ \ \ \ \ \ \ \ {\rm (see\
		(\ref{volume-comp-altalanos-0}))} \\
	&\leq&\limsup_{r\to \infty}\frac{{\rm Vol}_g(B(x,r))}{r^n} \ \ \ \ \
	\ \ \ \ \ \ \ \ \ \ \ \
	\ \ \ \ \ \ \ \ \ \ \ \ \ \ \ \ \ \
	\ \ \ \ \ \ \ {\rm (monotonicity)}\\&\leq& \limsup_{r\to
		\infty}\frac{{\rm Vol}_g(B(x_0,r+d({x_0},x)))}{r^n}\ \ \ \ \ \ \ \ \
	\ \ \ \  (B(x,r)\subset
	B(x_0,r+d({x_0},x))) \\
	&=& \limsup_{r\to
		\infty}\left(\frac{{\rm Vol}_g(B(x_0,r+d({x_0},x)))}{(r+d({x_0},x))^n}\cdot\frac{(r+d({x_0},x))^n}{r^n}\right) \\
	&=&\omega_n. \ \ \ \ \ \ \ \ \ \ \
	\ \ \ \ \ \ \ \ \ \ \ \ \ \ \ \ \ \
	\ \ \ \ \ \ \ \ \ \ \ \ \ \ \ \ \ \ \ \ \ \ \ \ \ \ \ \ \ \ \ \  \ \ \ \ \  \ {\rm
		(see\ (\ref{ideiskell}))}
\end{eqnarray*}
Consequently,
\begin{equation}\label{vol-id-neg}
{{\rm Vol}_g(B(x,\rho))}= \omega_n \rho^n\ {\rm for\ all}\ x\in M\ {\rm
	and}\ \rho>0.
\end{equation}
Now, the equality case in Theorem \ref{comparison-volume}(a) implies
that the sectional curvature is identically zero, which concludes
the proof. \hfill $\square$


\begin{remark}\rm
	Implication (a)$\Rightarrow$(c) in Theorem
	\ref{theorem-uncertainty-intro-1} has also a geometric
	proof.
	Indeed, due to Jost \cite[Lemma 2.1.5]{Jost} and relation
	(\ref{delta-egyenloseg}), it follows that we have equality in the
	CAT(0)-inequality with reference point $x_0\in M$, i.e., for {\it
		every} geodesic segment $\gamma:[0,1]\to M$ and $s\in [0,1]$, we
	have
	$$d^2(x_0,\gamma(s))= (1-s)d^2(x_0,\gamma(0))+sd^2(x_0,\gamma(1))-s(1-s)d^2(\gamma(0),\gamma(1)).$$
	Now, Alexandrov's rigidity result implies that the geodesic triangle
	formed by the points $x_0$, $\gamma(0)$ and $\gamma(1)$ is {\it
		flat}, see e.g. Bridson and Haefliger \cite{Br-Ha}. Therefore, the
	conclusion that $(M,g)$ is isometric to the Euclidean space $\mathbb
	R^n$ follows in a standard manner; the author thanks J. Jost and A.
	Lytchak for pointing out this approach.
\end{remark}


\subsection{Sharp Heisenberg-Pauli-Weyl principle on hyperbolic
	spaces}\label{subsec-4-3} For the hyperbolic space we use the
Poincar\'e ball model $\mathbb H^n=\{x\in \mathbb R^n:|x|<1\}$
endowed with the Riemannian metric $$g_{\rm
	hyp}(x)=(g_{ij}(x))_{i,j={1,...,n}}=p(x)^2\delta_{ij},$$ where
$p(x)=\frac{2}{1-|x|^2}.$ It is well known  that $(\mathbb
H^n,g_{\rm hyp})$ is a Cartan-Hadamard manifold with constant
sectional curvature $-1$. The volume form is
\begin{equation}\label{volume-form-hyper}
{\text d}V_{\mathbb H^n}(x) = p(x)^n {\text d}x,
\end{equation}
while the hyperbolic gradient and Laplace-Beltrami operator are
given by
$$\nabla_{\mathbb H^n}u=\frac{\nabla u}{p^2}\ {\rm and}\ \Delta_{\mathbb H^n}u=p^{-n}{\rm div}(p^{n-2}\nabla u),$$
where $\nabla$ denotes the Euclidean gradient in $\mathbb R^n.$ The
hyperbolic distance between the origin and $x\in \mathbb H^n$ is
given by $$d_{\mathbb
	H^n}(0,x)=\ln\left(\frac{1+|x|}{1-|x|}\right).$$

Recently, Kombe and \"Ozaydin \cite{KO-TAMS-2013} stated a
Heisenberg-Pauli-Weyl principle on $(\mathbb H^n,g_{\rm hyp})$. For
completeness, we recall the statement of Theorem 4.2 from \cite{KO-TAMS-2013}:\\

\noindent {\it $"$Let $u\in C_0^\infty(\mathbb H^n)$, $d=d(x)=d_{\mathbb
		H^n}(0,x)$ and $n>2$. Then}
\begin{equation}\label{kombe-up}
\left(\int_{\mathbb H^n}|\nabla_{\mathbb H^n} u|^2{\text d}V_{\mathbb H^n}\right) \left(\int_{\mathbb
	H^n}d^2u^2{\textrm d}V_{\mathbb H^n}\right)\geq
\frac{n^2}{4}\left(\int_{\mathbb
	H^n}u^2{\text d}V_{\mathbb H^n}\right)^2.
\end{equation}
{\it Moreover, equality holds in $(\ref{kombe-up})$ if
	$u(x)=Ae^{-\alpha d^2}$, where $A\in \mathbb R$, and
	\begin{equation}\label{rossz-egyenlet}
	\alpha=\frac{n-1}{n-2}\left(n-1+2\pi \frac{C_{n-2}}{C_n}\right)
	\end{equation}
	with} $C_n=\displaystyle\int_{\mathbb H^n}e^{-\alpha
	d^2}{\text d}V_{\mathbb H^n},$ $\alpha>0."$ \\

\noindent Relation (\ref{kombe-up}) holds true, see also Theorem \ref{theorem-uncertainty-intro-1}. However, the statement concerning the
{\it equality} in $(\ref{kombe-up})$ cannot happen, which has the following
three independent proofs:\\

\noindent 	{\bf Argument 1} (based on the non-solvability of (\ref{rossz-egyenlet}))	.
Let $C_n=C_n(\alpha)=\displaystyle\int_{\mathbb H^n}e^{-\alpha d^2}{\text d}V_{\mathbb
	H^n}$ be as above. We claim that the 
non-linear equation (\ref{rossz-egyenlet})  {\it cannot} be solved generically in $\alpha>0$. For simplicity, we consider only the case $n=4$; then equation (\ref{rossz-egyenlet}) reduces to $\alpha=w(\alpha)$, where 
$$w(\alpha):=\frac{3}{2}\left(3+2\pi
\frac{\displaystyle\int_{\mathbb H^2}e^{-\alpha d^2}{\text d}V_{\mathbb
		H^2}}{\displaystyle\int_{\mathbb H^4}e^{-\alpha d^2}{\text d}V_{\mathbb
		H^4}}\right).$$
Since $w\geq \frac{9}{2}$, the values for $\alpha$ should belong to $[\frac{9}{2},\infty)$ in order to solve $\alpha=w(\alpha)$. 

We claim that 
\begin{equation}\label{w-alpha}
w(\alpha)\geq 2\alpha+1\ \ {\rm  for\ every}\ \alpha\in [4,\infty),
\end{equation}
which will clearly imply the non-solvability of  $\alpha=w(\alpha)$. 

By (\ref{volume-form-hyper}), a change of variables shows that 
\begin{align*} w(\alpha)=\frac{9}{2}+\frac{\displaystyle 3\int_0^\infty e^{-\alpha t^2}\sinh(t){\rm d}t}{\displaystyle\int_0^\infty e^{-\alpha t^2}\sinh^3(t){\rm d}t}&=\frac{9}{2}+\frac{12 \erf\left(\frac{1}{2\sqrt{\alpha}}\right)}{\ds e^{\frac{2}{\alpha}}\erf\left(\frac{3}{2\sqrt{\alpha}}\right)-3\erf\left(\frac{1}{2\sqrt{\alpha}}\right)},
\end{align*}
where $\erf(s)=\ds \frac{2}{\sqrt{\pi}}\int_0^s e^{-t^2}{\rm d}t$ is the Gauss error function. 
Therefore, the claim (\ref{w-alpha}) is equivalent to the inequality $$3 \frac{4\alpha+1}{4\alpha-7}e^{-\frac{2}{\alpha}}{\ds \erf\left(\frac{1}{2\sqrt{\alpha}}\right)}\geq {\ds \erf\left(\frac{3}{2\sqrt{\alpha}}\right)}.$$ If $s=\ds \frac{1}{2\sqrt{\alpha}}\in \left(0,\frac{1}{4}\right]$, the latter inequality is equivalent to \begin{equation*}
\label{igazolni}
3\ds \frac{1+s^2}{1-7s^2} e^{-8s^2}\geq \ds\frac{\erf(3s)}{\erf(s)},\quad s\in \left(0,\frac{1}{4}\right].
\end{equation*}
Simple estimates for the error and exponential functions give for every $s\in \left(0,\frac{1}{4}\right]$ that $$3\frac{1+s^2}{1-7s^2}\cdot e^{-8s^2}-\frac{\erf(3s)}{\erf(s)}\geq 3{\frac {\ds \left( 1+{s}^{2} \right)  \left(1-4\,{s}^{2}\right)^2 }{1-
		7\,{s}^{2}}}-3{\frac {1-\,{s}^{2}\,}{1-\frac{1}{3}\,{s}^{2}}}\geq 0,
$$
which concludes the proof of  (\ref{w-alpha}). \\

\noindent 	{\bf Argument 2} (based on Theorem \ref{theorem-uncertainty-intro-1}). Following Kombe and \"Ozaydin \cite{KO-TAMS-2013}, let us assume that the hyperbolic Gaussian $u=e^{-\alpha d^2}>0$ is an extremal
function in $(\ref{kombe-up})$ for some $\alpha>0$. Due to Theorem
\ref{theorem-uncertainty-intro-1} (ii), it follows that   the hyperbolic space $(\mathbb H^n,g_{\rm hyp})$ is isometric to the standard Euclidean space $\mathbb R^n$, a contradiction.\\

\noindent 	{\bf Argument 3} (based on Theorem \ref{quantitative-HPW}).	  Due to Theorem \ref{quantitative-HPW}, for every $u\in C_0^\infty(\mathbb H^n)$ one has
\begin{equation}\label{uncer-qUANT1}
\left(\int_{\mathbb H^n}|\nabla_{\mathbb H^n} u|^2{\text d}V_{\mathbb H^n}\right) \left(\int_{\mathbb
	H^n}d^2u^2{\text d}V_{\mathbb H^n}\right)\geq
\frac{n^2}{4}\left(\int_{\mathbb
	H^n}\left(1+\frac{n-1}{n}{\bf D}_{-1}(d)\right)u^2{\text d}V_{\mathbb H^n}\right)^2.
\end{equation}
Since ${\bf D}_{-1}(d)\geq 0,$ if one expects to have equality in
(\ref{kombe-up}) for $u=e^{-\alpha d^2}$ for some $\alpha>0$, we necessarily have
in (\ref{uncer-qUANT1}) the relation ${\bf D}_{-1}(\rho)=0$ for
every $\rho\geq 0$; this relation means that for every $\rho\geq 0$ we have 
$$0=\rho {\bf ct}_{-1}(\rho)-1=\rho \coth (\rho)-1,$$
a contradiction.	Moreover, in the inequality (\ref{uncer-qUANT1}) the constant
$\frac{n^2}{4}$ is {\it sharp} and an integration by parts easily
shows (by using the exact form of the volume element
(\ref{volume-form-hyper})) that the {\it equality} holds for the
hyperbolic Gaussian family of functions $u_\alpha=e^{-\alpha
	d^2}$, $\alpha>0$. \\

Summing up the above discussions, we conclude that:
\begin{quote}
	{\it The  hyperbolic Gaussian functions $u_\lambda=e^{-\lambda
			d^2}$, $\lambda>0$,
		represent the family of extremals for the quantitative Heisenberg-Pauli-Weyl
		principle {\rm (\ref{uncer-qUANT1})}, but {\it not} for the 'pure'
		Heisenberg-Pauli-Weyl principle {\rm (\ref{kombe-up})}}.
\end{quote}


\vspace{0.2cm}

\subsection{Non-negatively curved case:  proof of Theorem \ref{theorem-uncertainty-intro-2}}

Implications 
(c)$\Rightarrow$(b)$\Rightarrow$(a) trivially hold. The proof of the implication (a)$\Rightarrow$(c) is divided into four steps. Let $x_0\in M$ be fixed. 

{\it \underline{Step 1}.} If $(M,g)$ is isometric to $\mathbb R^n$,
then  ${\bf ({HPW})}_{x_0}$  can be transformed into the inequality
(\ref{euklideszi-HPW}) for which the standard class of Gaussian
functions are extremals.

For later use, if we consider the function
$T:(0,\infty)\to \mathbb R$ defined by
$$T(\lambda)=\int_{\mathbb R^n} e^{-2\lambda |x|^2}{\text d}x,\ \lambda>0,$$
the equality for the family of extremals in (\ref{euklideszi-HPW})
can be rewritten to the form
\begin{equation}\label{uncertainty-identity-1}
-\lambda T'(\lambda)=\frac{n}{2} T(\lambda),\ \lambda>0.
\end{equation}
Moreover, by the layer cake representation and changing a variable,
one has the following representations which are used later:
\begin{equation}\label{hat-ez-eleg}
T(\lambda)=4\lambda\omega_n \int_0^\infty  \rho^{n+1} e^{-2\lambda \rho^2}
{\text d}\rho=\frac{2}{(2\lambda)^\frac{n}{2}} \omega_n
\int_0^\infty t^{n+1}e^{-t^2}{\text d}t.
\end{equation}

{\it \underline{Step 2}.} Let $x_0\in M$ be fixed. By our hypothesis,  ${\bf({HPW})}_{x_0}$ holds; in particular, $(M,g)$ cannot
be compact. We consider the class of functions
$$\tilde u_\lambda(x)=e^{-\lambda d_{x_0}(x)^2},\ \lambda>0.$$
Clearly,  the function $\tilde u_\lambda$ can be approximated by
elements from $C_0^\infty(M)$ for every $\lambda>0$. By inserting
$\tilde u_\lambda$ into ${\bf({HPW})}_{x_0}$, and using
$|{\nabla}_gd_{x_0}|=1$ a.e. on $M,$ we obtain that
\begin{equation}\label{uncertainty-ineq}
2\lambda\int_{M}d_{x_0}^2e^{-2\lambda d_{x_0}^2}{\text d}V_g\geq
\frac{n}{2}\int_{M}e^{-2\lambda{d_{x_0}^2}} {\text d}V_g,\
\lambda>0.
\end{equation}
We introduce the function
$\mathscr T:(0,\infty)\to \mathbb R$ defined by
$$\mathscr T(\lambda)=\int_{M}e^{-2\lambda{d_{x_0}^2}} {\text d}V_g,\ \lambda>0.$$
By the layer cake representation,  $\mathscr T$ can be equivalently
rewritten to
\begin{eqnarray}
\mathscr T(\lambda) &=& \int_0^\infty {\rm
	Vol}_g\left(\left\{x\in M:e^{-{2\lambda} d_{x_0}^2 }>t\right\}\right){\text
	d}t=\int_0^1 {\rm Vol}_g\left(\left\{x\in M:e^{-{2\lambda} d_{x_0}^2 }>t\right\}\right){\text d}t \nonumber\\ &=&
4\lambda\int_0^\infty  {\rm
	Vol}_g(B(x_0,\rho))\rho e^{-2\lambda \rho^2} {\text d}\rho.\nonumber
\end{eqnarray}
Since the Ricci curvature is non-negative, one account of
(\ref{volume-comp-altalanos-2}), the function $\mathscr T$ is well
defined and differentiable. Thus, relation (\ref{uncertainty-ineq})
is equivalent to
\begin{equation}\label{T-ben-egyenlotlenseg}
-\lambda \mathscr T'(\lambda)\geq \frac{n}{2} \mathscr T(\lambda),\ \lambda>0.
\end{equation}

{\it \underline{Step 3}.} We shall prove that
\begin{equation}\label{step-67}
\mathscr T (\lambda)\geq  T(\lambda) \ {\rm for\ all}\ \lambda>0.
\end{equation}
By (\ref{uncertainty-identity-1}) and (\ref{T-ben-egyenlotlenseg})
it turns out that
$$\frac{\mathscr T'(\lambda)}{\mathscr T(\lambda)} \leq \frac{
	T'(\lambda)}{T(\lambda)},\ \lambda>0.$$ Integrating this inequality,
it yields that the function $\lambda\mapsto \frac{\mathscr
	T(\lambda)}{T(\lambda)}$ is non-increasing; in particular, for every
$\lambda>0,$
\begin{equation}\label{inf-infty-tk}
\frac{\mathscr
	T(\lambda)}{T(\lambda)}\geq \liminf_{\lambda\to \infty}
\frac{\mathscr T(\lambda)}{T(\lambda)}.
\end{equation}
Now, we shall prove  that
\begin{equation}\label{egy99-1}
\liminf_{\lambda\to
	\infty}\frac{\mathscr T(\lambda)}{T(\lambda)}\geq1.
\end{equation}
Due to relation (\ref{volume-comp-nullaban}), for every
$\varepsilon>0$ one can find $\rho_\varepsilon>0$ such that
$${\rm Vol}_g(B(x_0,\rho))\geq (1-\varepsilon)\omega_n\rho^n \ {\rm
	for\ all}\ \rho\in [0,\rho_\varepsilon].$$ Consequently, one has
\begin{eqnarray*}
	\mathscr T(\lambda) &=&
	4\lambda\int_0^\infty  {\rm Vol}_g(B(x_0,\rho))\rho e^{-2\lambda \rho^2} {\text d}\rho\ \ \ \ \  \ \  \ \ \ \  \ \  \ \ \ \  \ \ \  \\
	&\geq&  4\lambda(1-\varepsilon)\omega_n\int_0^{\rho_\varepsilon} \rho^{n+1}e^{-2\lambda \rho^2}{\text d}\rho \\
	&=&\frac{2}{(2\lambda)^\frac{n}{2}}
	(1-\varepsilon)\omega_n\int_0^{\sqrt{2\lambda}\rho_\varepsilon} t^{n+1}e^{-t^2}{\text d}t.\ \ \ \ \  \ \  \ \ \  \ \  \ \ \ \  \ \  (\sqrt{2\lambda}\rho=t)
\end{eqnarray*}
Now, by (\ref{hat-ez-eleg}), it yields that $$
\liminf_{\lambda\to
	\infty}\frac{\mathscr T(\lambda)}{T(\lambda)}\geq 1-\varepsilon. $$
Since $\varepsilon>0$ is arbitrary, relation
(\ref{egy99-1}) holds, so (\ref{inf-infty-tk}). This ends the proof
of the claim (\ref{step-67}).

{\it \underline{Step 4}.} Via (\ref{hat-ez-eleg}) and the
representation of $\mathcal T$, relation (\ref{step-67}) is
equivalent to
$$\int_0^\infty \left({\rm
	Vol}_g(B(x_0,\rho))-\omega_n\rho^n\right)\rho e^{-2\lambda
	\rho^2}{\text d}\rho\geq 0\ {\rm for\ all}\  \lambda>0.$$ Due to
(\ref{volume-comp-altalanos-2}), we have
\begin{equation}\label{vol-id-poz-uj}
{{\rm Vol}_g(B(x_0,\rho))}= \omega_n \rho^n\ {\rm for\ all}\ \ \rho>0.
\end{equation}

Now, let $x\in M$ and $\rho>0$ be arbitrarily fixed. Note that by
Theorem \ref{comparison-volume}(b)  the function $r\mapsto
\frac{{\rm Vol}_g(B(x,r))}{r^n}$ is non-increasing on $(0,\infty)$.
Therefore, we have
\begin{eqnarray*}
	\omega_n &\geq &\frac{{\rm Vol}_g(B(x,\rho))}{\rho^n}\ \ \ \ \ \ \ \
	\ \ \
	\ \ \ \ \ \ \ \ \ \ \ \ \ \ \ \ \ \
	\ \ \ \ \ \ \ \ \ \ \ \ \ \ \ \ \ \ \ \ \ \ \ {\rm (see\
		(\ref{volume-comp-altalanos-2}))} \\
	&\geq&\limsup_{r\to \infty}\frac{{\rm Vol}_g(B(x,r))}{r^n} \ \ \ \ \
	\ \ \ \ \ \ \ \ \ \ \ \
	\ \ \ \ \ \ \ \ \ \ \ \ \ \ \ \ \ \
	\ \ \ \ \ \ \ {\rm (monotonicity)}\\&\geq& \limsup_{r\to
		\infty}\frac{{\rm Vol}_g(B(x_0,r-d(x_0,x)))}{r^n}\ \ \ \ \ \ \ \ \ \
	\ \ \  (B(x,r)\supset
	B(x_0,r-d(x_0,x))) \\
	&=& \limsup_{r\to
		\infty}\left(\frac{{\rm Vol}_g(B(x_0,r-d(x_0,x)))}{(r-d(x_0,x))^n}\cdot\frac{(r-d(x_0,x))^n}{r^n}\right) \\
	&=&\omega_n. \ \ \ \ \ \ \ \ \ \ \
	\ \ \ \ \ \ \ \ \ \ \ \ \ \ \ \ \ \
	\ \ \ \ \ \ \ \ \ \ \ \ \ \ \ \ \ \ \ \ \ \ \ \ \ \ \ \ \ \ \ \  \ \ \ \ \  \ {\rm
		(see\ (\ref{vol-id-poz-uj}))}
\end{eqnarray*}
Consequently, one has
\begin{equation}\label{javoj-jonnn}
{\rm Vol}_g(B(x,\rho))= \omega_n\rho^n\ {\rm for\ all}\ x\in M,\
\rho \geq 0.
\end{equation}
Thus, the equality case in Theorem \ref{comparison-volume}(b)
implies that the sectional curvature identically vanishes, which
conludes the proof.
\hfill $\square$

\vspace{0.5cm}

\section{Inequalities related  to the Heisenberg-Pauli-Weyl principle on Cartan-Hadamard
	manifolds}\label{sect-4}

\subsection{Caffarelli-Kohn-Nirenberg interpolation
	inequality:  proof of Theorem \ref{theorem-interpolation-intro}}

The proof is similar to Theorem \ref{theorem-uncertainty-intro-1}.

(i) Let $x_0\in M$ and $u\in C_0^\infty(M)$.  By Theorem
\ref{comparison-laplace}, we have $d_{x_0}{\Delta_g}d_{x_0}\geq
n-1.$ Consequently,
\begin{eqnarray}\label{ezt-is-kell}
\int_{M} \frac{|u|^p}{d_{x_0}^q}{\text d}V_g &\leq& \frac{1}{n-1} \int_{M} \frac{|u|^p}{d_{x_0}^{q-1}}{\Delta_g} d_{x_0} {\text d}V_g\\
&=& -\frac{1}{n-1} \int_{M} \left\langle\nabla_g\left(\frac{|u|^p}{d_{x_0}^{q-1}}\right),{\nabla}_g d_{x_0}\right\rangle {\text d}V_g \nonumber \\
&=& -\frac{p}{n-1} \int_{M} \frac{|u|^{p-2}u}{d_{x_0}^{q-1}}\left\langle\nabla_g|u|,{\nabla}_g d_{x_0}\right\rangle
{\text d}V_g \nonumber\\&& +\frac{q-1}{n-1} \int_{M} \frac{|u|^p}{d_{x_0}^{q}} |{\nabla}_g
d_{x_0}|^2
{\text d}V_g. \nonumber
\end{eqnarray}
Since $|{\nabla}_g
d_{x_0}|=1$, a reorganization of the above estimate
implies that $$
\frac{n-q}{p} \int_{M}
\frac{|u|^p}{d_{x_0}^{q}}{\text d}V_g\leq -\int_{M}
\frac{|u|^{p-2}u}{d_{x_0}^{q-1}}\left\langle\nabla_g|u|,{\nabla}_g
d_{x_0}\right\rangle
{\text d}V_g\leq \int_{M}
\frac{|u|^{p-2}u}{d_{x_0}^{q-1}}|\nabla_g u|
{\text d}V_g.
$$
By applying the Schwartz inequality, it yields the desired
inequality ${\bf(CKN)}_{x_0}$.

The proof of the sharpness of $\frac{(n-q)^2}{p^2}$ in
${\bf(CKN)}_{x_0}$ works in a similar manner as in Theorem
\ref{theorem-uncertainty-intro-1}, by exploiting the fact that in
the Euclidean setting the inequality ${\bf (CKN)}_{x_0}$ has the
form
$$
\left(\int_{\mathbb R^n}|\nabla u|^2{\text d}x\right) \left(\int_{\mathbb R^n}\frac{|u|^{2p-2}}{|x|^{2q-2}}{\text d}x\right) \geq \frac{(n-q)^2}{p^2}\left(\int_{\mathbb R^n}
\frac{|u|^p}{|x|^q} {\text d}x\right)^2, \ \forall u\in C_0^\infty(\mathbb
R^n),
$$
which has a class of positive extremals given in
(\ref{caff-extremal}).

(ii)  (a)$\Rightarrow$(c). According to the hypothesis,
$\frac{(n-q)^2}{p^2}$ is sharp and there exists a positive extremal
function $w_0$ in ${\bf (CKN)}_{x_0}$ for some $x_0\in M.$ In
particular, in relation (\ref{ezt-is-kell}) we should have the
equality
\begin{equation}\label{egyenlo-u_0-elso}
\int_{M} \frac{w_0^p}{d_{x_0}^q}{\text d}V_g = \frac{1}{n-1}
\int_{M} \frac{w_0^p}{d_{x_0}^{q-1}}{\Delta_g} d_{x_0} {\text d}V_g.
\end{equation}
Since $w_0>0$ and $d_{x_0}{\Delta}_gd_{x_0}\geq n-1$,
relation (\ref{egyenlo-u_0-elso}) implies that we necessarily have
$d_{x_0}{\Delta}_gd_{x_0}=n-1,$ thus $ {\Delta}_g(d_{x_0}^2)=2n.
$
The rest of the proof is similar to that of Theorem
\ref{theorem-uncertainty-intro-1}.
\hfill $\square$


\subsection{Hardy-Poincar\'e inequality: proof of Theorem \ref{theorem-Hardy-intro}}

Before to prove Theorem \ref{theorem-Hardy-intro}, we present a
quantitative version of the Hardy-Poincar\'e inequality on
Cartan-Hadamard manifolds.

\begin{theorem}\label{theorem-Hardy-intro-1}{\rm [Quantitative Hardy-Poincar\'e inequality]} Let $(M,g)$ be
	an  $n-$di\-men\-sional $(n\geq 3)$ Cartan-Hadamard manifold with
	sectional curvature bounded from above by $c\leq 0.$
	Then for every $x_0\in M$ and $u\in C_0^\infty(M)$ we have
	{ $$
		\int_{M}|\nabla_gu|^2\text{\emph{d}}V_g \geq   \frac{(n-2)^2}{4}\int_{M}\left(1+\frac{2(n-1)}{n-2}{\bf
			D}_c(d_{x_0})\right)\frac{u^2}{d_{x_0}^2}\text{\emph{d}}V_g.
		\eqno{{\bf
				({HP})}_{x_0}}
		$$}
	In addition, the constant $\frac{(n-2)^2}{4}$ is sharp and never
	achieved.
\end{theorem}

{\it Proof.} Let $x_0\in M$ and $u\in C_0^\infty(M)$ be
arbitrarily and fix $\gamma=\frac{n-2}{2}>0$. We consider the
function $v=d_{x_0}^\gamma u$. Thus, for  $u=d_{x_0}^{-\gamma} v$
one has
$$\nabla_g u=-\gamma d_{x_0}^{-\gamma-1}v
\nabla_gd_{x_0}+d_{x_0}^{-\gamma}\nabla_gv.$$ Therefore, it yields
\begin{eqnarray*}
	|\nabla_g u|^2
	&\geq& \gamma^2 d_{x_0}^{-2\gamma-2}v^2 |\nabla_g d_{x_0}|^2 -2\gamma d_{x_0}^{-2\gamma-1}v \langle \nabla_g d_{x_0},\nabla_g v \rangle.
\end{eqnarray*}
Since $|\nabla_g d_{x_0}|=1$ a.e. on $M,$ after integrating the latter
inequality, we obtain
\begin{equation}\label{hardy-nulla-marad-kesobb}
\int_{M}|\nabla_g u|^2{\text d}V_g  \geq
\gamma^2\int_{M}d_{x_0}^{-2\gamma-2}v^2{\text d}V_g+R_0,
\end{equation}
where
\begin{eqnarray*}
	R_0 &=& -2\gamma\int_M
	d_{x_0}^{-2\gamma-1}v \langle \nabla_g d_{x_0},\nabla_g v
	\rangle{\text d}V_g
	=\frac{1}{2} \int_M \langle \nabla_g(v^2), \nabla_g(d_{x_0}^{-2\gamma})\rangle {\text d}V_g\\
	&=&-\frac{1}{2} \int_M
	v^2{\Delta}_g(d_{x_0}^{-2\gamma}){\text d}V_g \\
	&=& \gamma \int_M
	v^2d_{x_0}^{-2\gamma-2}\left(-2\gamma-1+d_{x_0}{\Delta}_gd_{x_0}\right){\text d}V_g\\
	&\geq & \frac{(n-1)(n-2)}{2}\int_M \left(d_{x_0}{\bf
		ct}_c(d_{x_0})-1\right)\frac{u(x)^2}{d_{x_0}^2} {\text d}V_g, \ \ \
	\ \ \ \ {\rm (see\ Theorem\ \ref{comparison-laplace})}\\
	&= & \frac{(n-1)(n-2)}{2}\int_M {\bf
		D}_c(d_{x_0})\frac{u(x)^2}{d_{x_0}^2} {\text d}V_g,
\end{eqnarray*}
which completes the first part of the proof.

We shall prove in the sequel that $\gamma^2=\frac{(n-2)^2}{4}$ is
sharp in ${{\bf ({HP})}_{x_0}}$, i.e.,
\begin{equation}\label{hardy-sharp}
\frac{(n-2)^2}{4}=\inf_{u\in
	C_0^\infty(M)\setminus\{0\}}\frac{\displaystyle\int_{M}|\nabla_g
	u|^2{\text
		d}V_g}{\displaystyle\int_{M}\left(1+\frac{2(n-1)}{n-2}{\bf
		D}_c(d_{x_0})\right)\frac{ u^2}{d_{x_0}^2}{\text d}V_g}.
\end{equation}
Fix the numbers $R>r>0$ and a smooth cutoff function $\psi:M\to
[0,1]$ with supp$(\psi)=\overline{B(x_0,R)}$ and $\psi(x)=1$ for
$x\in B(x_0,r)$, and for every $\varepsilon>0,$ let
\begin{equation}\label{u-eps-mindenki-ezt}
u_\varepsilon=(\max\{\varepsilon,d_{x_0}\})^{-\gamma}.
\end{equation}
On one hand,
\begin{eqnarray*}
	I_1(\varepsilon)&:=&\int_{M}|\nabla_g(\psi u_\varepsilon)|^2{\text
		d}V_g\\ &=& \int_{B(x_0,r)}|\nabla_g(\psi u_\varepsilon)|^2{\text
		d}V_g + \int_{B(x_0,R)\setminus B(x_0,r)}|\nabla_g(\psi u_\varepsilon)|^2{\text d}V_g\\
	&=& \gamma^2\int_{B(x_0,r)\setminus  B(x_0,\varepsilon)}d_{x_0}^{-2\gamma-2}{\text
		d}V_g+\tilde I_1(\varepsilon),
\end{eqnarray*}
where the quantity
$$ \tilde I_1(\varepsilon)=\int_{B(x_0,R)\setminus B(x_0,r)}|\nabla_g(\psi u_\varepsilon)|^2{\text d}V_g$$ is finite and does {\it not} depend on
$\varepsilon>0$ whenever $\varepsilon<r.$ On the other hand,
\begin{eqnarray*}
	I_2(\varepsilon)&:=&\int_{M}\left(1+\frac{2(n-1)}{n-2}{\bf
		D}_c(d_{x_0})\right)\frac{(\psi u_\varepsilon)^2}{d_{x_0}^2}{\text d}V_g\\
	&\geq& \int_{M}\frac{(\psi u_\varepsilon)^2}{d_{x_0}^2}{\text d}V_g\\
	&\geq& \int_{B(x_0,r)\setminus  B(x_0,\varepsilon)}d_{x_0}^{-2\gamma-2}{\text
		d}V_g=:\tilde I_2(\varepsilon).
\end{eqnarray*}
By applying the layer cake representation, we deduce that for
$0<\varepsilon<r$, one has
\begin{eqnarray*}
	\tilde I_2(\varepsilon)&=&\int_{B(x_0,r)\setminus
		B(x_0,\varepsilon)}d_{x_0}^{-2\gamma-2}{\text d}V_g  =
	\int_{B(x_0,r)\setminus
		B(x_0,\varepsilon)}d_{x_0}^{-n}{\text d}V_g \\
	&\geq& \int_{r^{-n}}^{\varepsilon^{-n}}{\rm Vol}_g(B(x_0,\rho^{-\frac{1}{n}})){\text d}\rho \\
	&\geq& \omega_n\int_{r^{-n}}^{\varepsilon^{-n}}\rho^{-1}{\text
		d}\rho\ \ \ \ \ \ \ \ \ \ \ \ \ \ \ \ \ \ \ \ \ \ \ \ \ \ \ \ \ \ \
	\ \ \ \ \ \ \ \ \ \ \ \ \ \ \ \ \ \ \ \ \ \ \ \ \ {\rm (see\
		(\ref{volume-comp-altalanos-0}))}\\
	&=&n\omega_n(\ln r-\ln \varepsilon).
\end{eqnarray*}
In particular, $\lim_{\varepsilon\to 0^+}\tilde
I_2(\varepsilon)=+\infty.$ Thus, from the above relations it follows
that
\begin{eqnarray*}
	\frac{(n-2)^2}{4} &\leq & \inf_{u\in
		C_0^\infty(M)\setminus\{0\}}\frac{\displaystyle\int_{M}|\nabla_g
		u|^2{\text
			d}V_g}{\displaystyle\int_{M}\left(1+\frac{2(n-1)}{n-2}{\bf
			D}_c(d_{x_0})\right)\frac{ u^2}{d_{x_0}^2}{\text
			d}V_g} \\
	&\leq& \lim_{\varepsilon\to
		0^+}\frac{I_1(\varepsilon)}{I_2(\varepsilon)}
	\leq\lim_{\varepsilon\to
		0^+}\frac{\gamma^2\tilde I_2(\varepsilon)+\tilde
		I_1(\varepsilon)}{\tilde
		I_2(\varepsilon)}\\&=&\gamma^2=\frac{(n-2)^2}{4},
\end{eqnarray*}
which concludes the proof of (\ref{hardy-sharp}).

If we assume the function  $u_0\neq 0$ is an extremal in ${{\bf
		({HP})}_{x_0}}$, on one hand, due to
(\ref{hardy-nulla-marad-kesobb}) we have
\begin{equation}\label{utolsok-kozotti}
\int_M d_{x_0}^{-2\gamma}|\nabla_gv_0|^2{\text d}V_g=0,
\end{equation}
where $v_0=d_{x_0}^\gamma u_0$. By (\ref{utolsok-kozotti}) it
follows that $v_0$ is a constant function, thus
$u_0=c_0d_{x_0}^{-\gamma}$ for some $c_0\in \mathbb R\setminus
\{0\}.$ On the other hand, similar estimates as above show that
$$ \int_M |\nabla_gu_0|^2{\text d}V_g=\gamma^2\int_M \frac{u_0^2}{d_{x_0}^2}{\text d}V_g=c_0^2\gamma^2\int_M d_{x_0}^{-n}{\text d}V_g=+\infty,$$
i.e., $u_0\notin H^1(M,{\text d}V_g)$ and $\frac{u_0}{d_{x_0}}\notin
L^2(M,{\text d}V_g)$, a contradiction.
\hfill $\square$\\

{\it {Proof of Theorem \ref{theorem-Hardy-intro}}.}
By the continued
fraction representation of the function $\rho\mapsto \coth(\rho)$,
one has
$$\rho\coth(\rho)-1\geq \frac{3\rho^2}{\pi^2+\rho^2}\ \ {\rm for\ all}\ 
\rho>0.$$ Now,  the inequality follows at once from this estimate
and Theorem \ref{theorem-Hardy-intro-1}.\hfill $\square$

\begin{remark} \rm (i) Our arguments work also for {\it weighted}
	Hardy-Poincar\'e inequalities; for simplicity, we presented ${{\bf
			({HP})}_{x_0}}$ in its simplest form.
	
	(ii) Kombe and \"Ozaydin \cite{KO-TAMS-2009, KO-TAMS-2013}
	investigated the sharp constant in the Hardy-Poincar\'e inequality
	on the hyperbolic space $\mathbb H^n,$ $n\geq 3.$ As expected, they
	claimed that
	\begin{equation}\label{hardy-vajon-joe}
	\frac{(n-2)^2}{4}=\inf_{u\in C_0^\infty(\mathbb H^n)\setminus
		\{0\}}{\small \frac{\displaystyle\int_{\mathbb H^n}|\nabla_{\mathbb
				H^n} u|^2{\text d}V_{\mathbb H^n}}{\displaystyle\int_{\mathbb
				H^n}\frac{u^2}{d^2}{\text d}V_{\mathbb H^n}},}
	\end{equation}
	where the
	notations come from \S \ref{subsec-4-3}.  In order to prove
	(\ref{hardy-vajon-joe}), the authors used as test functions {\it
		only} those from  (\ref{u-eps-mindenki-ezt}) without coupling with
	an appropriate cutoff function (as in the proof of Theorem
	\ref{theorem-Hardy-intro-1}). Although the functions $u_\varepsilon$
	can be approximated by elements from $C_0^\infty(\mathbb H^n)$, the
	gap in \cite{KO-TAMS-2009, KO-TAMS-2013} appears due to the fact
	that $u_\varepsilon\notin H^1(\mathbb H^n,{\text d}V_{\mathbb H^n})$
	and $\frac{u_\varepsilon}{d}\notin L^2(\mathbb H^n,{\text
		d}V_{\mathbb H^n})$,  $\varepsilon>0.$ Indeed, simple computations
	show that for every $\varepsilon>0,$
	\begin{eqnarray*}
		\int_{\mathbb H^n}|\nabla_{\mathbb H^n} u_\varepsilon|^2{\text d}V_{\mathbb H^n}  &=& \left(\gamma+\varepsilon\right)^2\int_{\mathbb H^n\setminus B(0,1)}
		d^{-2\gamma-2\varepsilon-2}{\text d}V_{\mathbb H^n}\\&=&
		\left(\gamma+\varepsilon\right)^2n\omega_n\int_1^\infty
		t^{-n-2\varepsilon}(\sinh t)^{n-1}{\text d}t=+\infty,
	\end{eqnarray*}
	and
	\begin{eqnarray*}
		\int_{\mathbb H^n}\frac{u_\varepsilon^2}{d^2}{\text d}V_{\mathbb
			H^n}\geq \int_{\mathbb H^n\setminus B(0,1)}
		d^{-2\gamma-2\varepsilon-2}{\text d}V_{\mathbb
			H^n}=n\omega_n\int_1^\infty t^{-n-2\varepsilon}(\sinh
		t)^{n-1}{\text d}t=+\infty.
	\end{eqnarray*}
	
	(iii) Similar observation as in (ii) has been already made in  Yang,
	Su and Kong \cite{YSK-Kinaiak}. In \cite{YSK-Kinaiak}, the authors
	proved sharp Hardy and Rellich inequalities on Riemannian manifolds
	with negative sectional curvature. The novelty of our results
	(Theorem \ref{theorem-Hardy-intro} \& \ref{theorem-Hardy-intro-1})
	is that improvements appear quantitatively in terms of the
	sectional curvature.
\end{remark}

A similar argument as in Theorem \ref{theorem-Hardy-intro-1} leads
to the following improvement.

\begin{theorem}\label{theorem-Hardy-intro-2}{\rm [Double improved Hardy-Poincar\'e inequality]} Let $\Omega$ be a bounded open domain with smooth boundary in an
	$n-$di\-men\-sional $(n\geq 3)$ Cartan-Hadamard manifold with
	sectional curvature bounded from above by $c\leq 0.$
	If $x_0\in \Omega$ and $R>\sup_{x\in \Omega} d(x,x_0)$, then
	for all $u\in C_0^\infty(\Omega),$ { $$
		\int_{\Omega}|\nabla_gu|^2\text{\emph{d}}V_g \geq   \frac{(n-2)^2}{4}\int_{\Omega}\left(1+\frac{2(n-1)}{n-2}{\bf
			D}_c(d_{x_0})\right)\frac{u^2}{d_{x_0}^2}\text{\emph{d}}V_g
		+\frac{1}{4}R_\Omega, $$} where
	$$
	R_\Omega=\int_{\Omega}\left(1+{2(n-1)}\ln\left(\frac{eR}{d_{x_0}}\right){\bf
		D}_c(d_{x_0})\right)\frac{u(x)^2}{d_{x_0}^2\ln^2\left(\frac{eR}{d_{x_0}}\right)}\text{\emph{d}}V_g.
	$$
\end{theorem}

 {\it Proof.} Let $x_0\in \Omega$, $u\in C_0^\infty(\Omega)$ and fix
$\gamma=\frac{n-2}{2}>0$. If we consider the function
$v=d_{x_0}^\gamma u$, one has
\begin{eqnarray*}
	|\nabla_gu|^2
	&=& \gamma^2 d_{x_0}^{-2\gamma-2}v^2-2\gamma  d_{x_0}^{-2\gamma-1}v\langle \nabla_g d_{x_0},\nabla_g v \rangle
	+d_{x_0}^{-2\gamma}|\nabla_gu|^2.
\end{eqnarray*}
After an integration over $\Omega$ of the above relation, one can
repeat the argument from the proof of Theorem
\ref{theorem-Hardy-intro-1} to the first two integrands, obtaining
$$\int_{\Omega}|\nabla_gu|^2{\text d}V_g  \geq  \frac{(n-2)^2}{4}\int_{\Omega}\left(1+\frac{2(n-1)}{n-2}{\bf
	D}_c(d_{x_0})\right)\frac{u^2}{d_{x_0}^2}{\text d}V_g+\tilde R,$$
where
$$\tilde R=\int_\Omega d_{x_0}^{-2\gamma}|\nabla_gv|^2{\text d}V_g.$$
Due to the fact that $R>\sup_{x\in \Omega} d(x,x_0)$, the function
$h=\ln \frac{eR}{d_{x_0}}$ is well defined on $\Omega\setminus
\{x_0\}$ and $h\geq 1.$ Let $z=h^{-1/2}v.$ Since
$$\nabla_gv=-\frac{z}{2d_{x_0}}h^{-1/2}\nabla_gd_{x_0}+h^{1/2}\nabla_gz,$$
it turns out that
$$|\nabla_gv|^2\geq \frac{z^2}{4d_{x_0}^2}h^{-1}-\frac{z}{d_{x_0}}\langle \nabla_gd_{x_0}, \nabla_gz\rangle.$$
Consequently,
\begin{eqnarray*}
	\tilde R &=& \int_\Omega d_{x_0}^{-2\gamma}|\nabla_gv|^2{\text d}V_g \\
	&\geq&\frac{1}{4} \int_\Omega d_{x_0}^{-2\gamma-2} h^{-1}z^2{\text d}V_g-\frac{1}{2}\int_\Omega d_{x_0}^{-2\gamma-1}\langle\nabla_gd_{x_0}, \nabla_g(z^2)\rangle{\text d}V_g\\
	&=&\frac{1}{4} \int_\Omega d_{x_0}^{-2} h^{-2}u^2{\text d}V_g-\frac{1}{4\gamma}\int_\Omega z^2{\Delta}_g(d_{x_0}^{-2\gamma}){\text d}V_g\\
	&=&\frac{1}{4} \int_\Omega d_{x_0}^{-2} h^{-2}u^2{\text d}V_g+\frac{1}{2}\int_\Omega z^2d_{x_0}^{-2\gamma-2}(-2\gamma-1+d_{x_0}{\Delta}_gd_{x_0}){\text d}V_g\\
	&\geq& \frac{1}{4} \int_\Omega d_{x_0}^{-2} h^{-2}u^2{\text
		d}V_g+\frac{n-1}{2}\int_\Omega {\bf
		D}_c(d_{x_0})d_{x_0}^{-2\gamma-2} z^2{\text d}V_g\\&=&
	\frac{1}{4}\int_{\Omega}\left(1+{2(n-1)}h{\bf
		D}_c(d_{x_0})\right)\frac{u^2}{d_{x_0}^2h^2}{\text d}V_g,
\end{eqnarray*}
which concludes the proof.  \hfill $\square$

\begin{remark}\rm  In the limiting case when $c=0$ (thus  ${\bf D}_c(\rho)={\bf
		D}_0(\rho)=0$ for every $\rho\geq 0$), the inequality in Theorem
	\ref{theorem-Hardy-intro-2} takes the familiar form
	$$\int_{\Omega}|\nabla_gu|^2\text{{d}}V_g \geq   \frac{(n-2)^2}{4}\int_{\Omega}\frac{u^2}{d_{x_0}^2}{\text d}V_g+\frac{1}{4}\int_{\Omega}\frac{u^2}{d_{x_0}^2\ln^2\left(\frac{eR}{d_{x_0}}\right)}{\text d}V_g,$$
	see Adimurthi,  Chaudhuri and Ramaswamy \cite{AR} and Filippas and
	Tertikas \cite{Fil-Ter} in the Euclidean case, and Kombe and
	\"Ozaydin \cite[Corollary 2.2]{KO-TAMS-2013} in hyperbolic spaces.
\end{remark}

\vspace*{0.2cm}
\noindent {\bf Acknowledgment.}  The author   would     like     to     thank     the   two  anonymous  reviewers for   their suggestions  and comments.


\end{document}